\documentclass{article}
\usepackage[T1]{fontenc} % check whether it works
\usepackage{relsize}
\usepackage[utf8]{inputenc}
\usepackage{realboxes}
\usepackage{dsfont}
\usepackage{soul}
\usepackage{marvosym}
\usepackage[normalem]{ulem}
\usepackage{leftindex}
\usepackage{float}
\usepackage{dashbox}
\usepackage{fancybox}
\usepackage[english]{babel}
\usepackage{graphicx}%
\usepackage{multirow}%
\usepackage{amsmath,amssymb,amsfonts}%
\usepackage{amsthm}%
\usepackage{mathrsfs}%
\usepackage[title]{appendix}%
\usepackage[svgnames]{xcolor}%
\usepackage{textcomp}%
\usepackage{manyfoot}%
\usepackage{booktabs}%
\usepackage{algorithm}%
\usepackage{algorithmicx}%
\usepackage{algpseudocode}%
\usepackage{listings}%
\usepackage{tikz-cd}
\usepackage{pgfplots}
\pgfplotsset{compat=1.13}
\usepackage[all,cmtip]{xy}
\usepackage{graphicx,tikz-cd,pgf}
\usetikzlibrary{positioning, shapes}
\usetikzlibrary{arrows.meta}
\usepackage{caption}
\usepackage{tkz-euclide}

\usepackage{stmaryrd} \usepackage{trimclip}
\makeatletter
\DeclareRobustCommand{\shortto}{\mathrel{\mathpalette\short@to\relax}}
\newcommand{\short@to}[2]{
  \mkern2mu
  \clipbox{{.5\width} 0 0 0}{$\m@th#1\vphantom{+}{\shortrightarrow}$}}
\makeatother

\makeatletter
\DeclareRobustCommand{\lshortto}{\mathrel{\mathpalette\lshort@to\relax}}
\newcommand{\lshort@to}[2]{
  \mkern2mu
  \clipbox{{.1\width} 0 0 0}{$\m@th#1\vphantom{+}{\shortleftarrow}$}}
\makeatother

\newcommand{\cfunction}[1]{{\mathtt{1}}_{#1}}

\DeclareMathOperator{\PoissonISymbol}{\mathsf{P}}
\newcommand{\PoissonI}[1]{\PoissonISymbol[#1]}

\DeclareMathOperator{\oscSymbol}{osc}
\newcommand{\osc}[2]{\oscSymbol(#1,#2)}

\DeclareMathOperator{\nullsets}{\boldsymbol{\mathcal{N}}}

\DeclareMathOperator{\almostev}{\text{\tiny{a.e.}}}

\newcommand{\soba}{\mathsf{Arcs}}

\newcommand{\bafNAME}{\boldsymbol{\Lambda}}
\newcommand{\baf}[1]{
\bafNAME(#1)
}
\newcommand{\bafz}{\bafNAME}

\newcommand{\bafpaz}{% PA, Zero variables
\bafNAME_{\pbsf}
}

\newcommand{\baff}[1]{% Firm
\bafNAME_{\fbsf}(#1)
}

\newcommand{\bafaz}{% Attached, Zero variables
\bafNAME_{\absf}
}

\newcommand{\bafcz}{% curvilinear, Zero variables
\bafNAME_{\cbsf}
}

\newcommand{\bafc}[1]{
\bafNAME_{\cbsf}(#1)
}

\newcommand{\bafsz}{% sequential, Zero variables
\bafNAME_{\sbsf}
}

\newcommand{\bafsNAME}{\bafNAME_{\sbsf}}
\newcommand{\bafs}[1]{
\bafsNAME(#1)
}

\DeclareMathOperator{\FsSymbol}{F}
\newcommand{\Fs}[1]{% Fatou Set
\FsSymbol[#1]
}

\newcommand{\baftNAME}{\bafNAME_{\tbsf}}
\newcommand{\baft}[1]{
\baftNAME(#1)
}

\newcommand{\bafoNAME}{\bafNAME_{\obsf}}
\newcommand{\bafo}[1]{
\bafoNAME(#1)
}

\newcommand{\baftz}{% Zero variables
\baftNAME
}

\DeclareMathOperator{\bafmNAME}{R}
\newcommand{\bafm}[1]{
\bafmNAME(#1)
}

\newcommand{\aesubset}[2]{
#1
\overset{
\scriptstyle{\almostev}
}
{\subset}
#2
}

\newcommand{\domain}{\Omega}

\newcommand{\bdryD}{\partial\domain}

\newcommand{\bdryhd}{\partial\hd}

\newcommand{\aeequal}[2]{
#1
\overset{
\scriptstyle{\almostev}
}
{=}
#2
}

\newcommand{\acal}[2]{% arc (with) center (and) arc-length
\mathsf{arc}[#1;#2]
}

\DeclareMathOperator{\AcalSymbol}{\mathtt{arcs}}

\DeclareMathOperator{\AcalSymbolS}{\mathtt{Arcs}}

\newcommand{\Acal}[2]{% Union (of ) arcs 
{\cup}_{\AcalSymbol}[#1;#2]
}

\newcommand{\racal}[2]{% right arc (with) center (and) arc-length
{\mathtt{arc}}_{+}[#1;#2]
}

\newcommand{\iweaw}{\AcalSymbolS{\!(\bpoint)}}

\newcommand{\iwleaw}{\AcalSymbolS\!{(\bpoint\!\!\lshortto\!)}}

\newcommand{\carcat}[1]{\mathsf{#1}}

\newcommand{\sCattop}{\boldsymbol{
{\carcat{top}}}}

\newcommand{\ichar}[1]{\mathtt{#1}}
\newcommand{\Ui}{\ichar{U}}

\newcommand{\Qb}{\bchar{Q}}
\newcommand{\Eb}{\bchar{E}}
\newcommand{\Nb}{\bchar{N}}
\newcommand{\yb}{\bchar{y}}

\newcommand{\Xb}{\bchar{X}}
\newcommand{\Yb}{\bchar{Y}}

\newcommand{\Jb}{\bchar{J}}
\newcommand{\Bb}{\bchar{B}}
\newcommand{\Vb}{\bchar{V}}
\newcommand{\Wb}{\bchar{W}}
\newcommand{\Cb}{\bchar{C}}

\newcommand{\Ab}{\bchar{A}}

\newcommand{\Ob}{\bchar{O}}
\newcommand{\Sb}{\bchar{S}}
\newcommand{\Ub}{\bchar{U}}

\newcommand{\vb}{\bchar{v}}
\newcommand{\Zb}{\bchar{Z}}

\newcommand{\lingvalue}[3]{% 1 function 2 approach region 3 boundary point
\lim_{\stackrel{z\in #2}{z\to #3}}{#1(z)}
}

\newcommand{\limitingv}[3]{\lim_{{#2{\shortto}#1}}#3}

\newcommand{\inter}{\between}

\DeclareMathOperator{\reg}{\boldsymbol{\mathsf{reg}}}

\DeclareMathOperator{\invariant}{\boldsymbol{\mathsf{invariant}}}

\DeclareMathOperator{\snY}{Y}

\DeclareMathOperator{\subbase}{\mathcal{C}}

\DeclareMathOperator{\osubbase}{\mathsf{C}}

\newcommand{\grigio}[1]{\colorbox{gray!20}
{
\ensuremath{\displaystyle{#1}}
}
}

\definecolor{mygray}{rgb}{0.8,0.8,0.8}
\makeatletter
    \newcommand{\colorboxed}[3][white]{\fcolorbox{#2}{#1}{\m@th$\displaystyle#3$}}
\makeatother

\DeclareMathOperator{\approachrS}{S}

\DeclareMathOperator{\dpoint}{\ichar{z}}

\DeclareMathOperator{\bpoint}{\mathsf{w}}
\DeclareMathOperator{\bopoint}{\bpoint}
\newcommand{\dsegnato}{{d}_{o}}
\newcommand{\bsegnato}{{b}^{\prime}}
\newcommand{\zsegnato}{\tilde{\dpoint}}
\DeclareMathOperator{\bpointv}{\mathsf{v}}
\DeclareMathOperator{\bpointx}{\mathsf{x}}

\DeclareMathOperator{\contact}{\boldsymbol{\tau}}

\newcommand{\utau}[1]{\contact^{*}\!\approachr(#1)}
\newcommand{\ltau}[1]{\contact_{*}\!\approachr(#1)}

\newcommand{\utauprime}[1]{\contact^{*}\!{\approachrprime}(#1)}
\newcommand{\ltauprime}[1]{\contact_{*}\!{\approachrprime}(#1)}

\DeclareMathOperator{\ballName}{{{B}}}
\newcommand{\ball}[2]{\ballName(#1,#2)}
\newcommand{\afballrel}[2]{\af(#1,#2)}
\newcommand{\ballrel}[2]{\approachr_{#1}^{[#2]}}
\newcommand{\ballrelprime}[2]{\approachrprime(#1,#2)}

\newcommand{\ballrelprimeprime}[2]{\approachrprimeprime(#1,#2)}

\DeclareMathOperator{\coronaName}{C}
\newcommand{\corona}[3]{\coronaName(#1;#2,#3)}

\DeclareMathOperator{\ntbvSymbol}{
\boldsymbol{\flat}
}
\newcommand{\ntbv}[1]{
{#1}_{\ntbvSymbol}
}

\newcommand{\angolo}{b}

\makeatletter
\newcommand*\bigcdot{\mathpalette\bigcdot@{.5}}
\newcommand*\bigcdot@[2]{\mathbin{\vcenter{\hbox{\scalebox{#2}{$\m@th#1\bullet$}}}}}
\makeatother

\DeclareMathSymbol{\mhyphen}{\mathord}{AMSa}{"39}

\DeclareMathOperator{\limite}{lim}

\DeclareMathOperator{\Limite}{{\mathtt{lim}}}
\newcommand{\Limiteop}[2]{\Limite_{#1}#2}

\DeclareMathOperator{\cl}{\mathsf{cl}}
\DeclareMathOperator{\clinfty}{\cl_{\infty}}
\newcommand{\clos}[1]{\cl(#1)}
\newcommand{\closinfty}[1]{\clinfty(#1)}

\DeclareMathOperator{\udone}{%BoundaryOfTheUnitDisc
\mathbb{D}
}

\DeclareMathOperator{\botud}{%BoundaryOfTheUnitDisc
\partial\udone
}

\DeclareMathOperator{\edome}{E}
\DeclareMathOperator{\edomd}{d}
\DeclareMathOperator{\edom}{{\edome}_{\edomd}}
\newcommand{\newedom}[1]{{\edom}(#1)}
\newcommand{\edomain}[1]{\edom(#1)}

\usepackage[mathscr]{euscript}

\usepackage{stmaryrd} \usepackage{trimclip}

\newcommand{\parrow}
{
{
\,\,\normalfont\scalebox{1.3}{$\rightarrowtail$}
\,
}
}

\DeclareMathOperator{\interval}{\Jb}

\newcommand{\res}{% residual
\rbsf[\bpoint]}

\newcommand{\areaw}{% a(pproach) r(egions) e(nding) a(t) w
\udone[\bpoint]}

\newcommand{\areawprime}{% a(pproach) r(egions) e(nding) a(t) w
\udone[\bpoint^{\prime}]}

\newcommand{\aretaw}{% 
\tbsf[\bpoint]}

\newcommand{\arenaw}{% 
{\Stolz}[\bpoint]}

\newcommand{\arvoaw}{% 
\obsf[\bpoint]}

\newcommand{\arpaw}{% 
\cbsf[\bpoint]}

\newcommand{\arcaw}{% 
\absf[\bpoint]}

\newcommand{\arpaaw}{%  
\pbsf[\bpoint]}

\newcommand{\varpaaw}[1]{%  
\pbsf_{#1}[\bpoint]}

\newcommand{\arsaw}{%  
\sbsf[\bpoint]}

\newcommand{\uniform}{\boldsymbol{\mathsf{unif}}}

\newcommand{\approach}{\boldsymbol{\ell}}
\newcommand{\Stolz}{\boldsymbol{\gamma}}
\newcommand{\tangential}{\tbsf}

\newcommand{\absf}{\boldsymbol{\mathsf{a}}}

\newcommand{\bbsf}{\boldsymbol{\mathsf{b}}}
\newcommand{\cbsf}{\boldsymbol{\mathsf{c}}}

\newcommand{\dbsf}{\boldsymbol{\mathsf{d}}}
\newcommand{\fbsf}{\boldsymbol{\mathsf{f}}}

\newcommand{\obsf}{\boldsymbol{\mathsf{o}}}

\newcommand{\pbsf}{\boldsymbol{\mathsf{p}}}

\newcommand{\rbsf}{\boldsymbol{\mathsf{r}}}
\newcommand{\sbsf}{\boldsymbol{\mathsf{s}}}

\newcommand{\tbsf}{\boldsymbol{\mathsf{t}}}

\newcommand{\vbsf}{\boldsymbol{\mathsf{v}}}
\newcommand{\xbsf}{\boldsymbol{\mathsf{x}}}

\newcommand{\zbsf}{\boldsymbol{\mathsf{z}}}

\newcommand{\spro}{% standard product
{\botud\times\tps{\udone}}
}

\newcommand{\sofor}{% s(et) o(f) f(amilies) o(f) r(egions)
\shs{\botud}{\tps{\udone}}}

\newcommand{\sofoar}{% s(et) o(f) f(amilies) o(f) a(pproach) r(egions)
{\sofor}_{\approach}
}

\newcommand{\sofotar}{% 
{\sofor}_{\tangential}
}

\newcommand{\sofoSar}{% 
{\sofor}_{\Stolz}
}

\newcommand{\sofovoar}{% 
{\sofor}_{\obsf}
}

\newcommand{\sofopar}{% 
{\sofor}_{\cbsf}
}

\newcommand{\sofocar}{% 
{\sofor}_{\absf}
}

\newcommand{\pt}{% 
{\sofor}_{\tbsf\cbsf}
}

\newcommand{\stnew}{
{\sofor}_{\tbsf\sbsf}
}

\newcommand{\op}{% 
{\sofor}_{\obsf\cbsf}
}

\newcommand{\sofoseqar}{% 
{\sofor}_{\sbsf}
}

\newcommand{\shs}[2]{
{
\boldsymbol{[\![}
#1
\shortto
#2
\boldsymbol{]\!]}
}
}

\newcommand{\setb}{B}

\newcommand{\setl}{L}

\newcommand{\continuous}[2]{
\cathomset{\sCattop}{#1}{#2}
}

\newcommand{\partialhomset}[2]{
{
\boldsymbol{\{\!\!\{}
#1
\shortto
#2
\boldsymbol{\}\!\!\}}
}
}

\newcommand{\cathomset}[3]{
{
\boldsymbol{[\![}
#2
\shortto
#3
\boldsymbol{]\!]}
}_{(#1)}
}

\newcommand{\dfunctionharm}{\mathtt{u}}
\newcommand{\dfunctionharmc}{\mathtt{v}}
\newcommand{\dfunctionhol}{\mathtt{h}}
\newcommand{\dfunction}{\mathtt{f}}
\newcommand{\bfunction}{\mathsf{f}}

\DeclareMathOperator{\approachr}{% approach r(egion)
\afn}
\DeclareMathOperator{\approachrprime}{% approach r(egion)
\mathsf{V}}
\DeclareMathOperator{\approachrprimeprime}{% approach r(egion)
\mathsf{W}}

\DeclareMathOperator{\af}{\boldsymbol{\mathsf{A}}}
\DeclareMathOperator{\afn}{ % n = normal 
{\mathsf{A}}}

\DeclareMathOperator{\afprime}{\boldsymbol{\mathsf{V}}}
\DeclareMathOperator{\afb}{\boldsymbol{\mathsf{W}}}

\DeclareMathOperator{\property}{\textsc{Q}}

\DeclareMathOperator{\regular}{R(\botud,\tps{\udone})}

\DeclareMathOperator{\tpsSymbol}{\mathcal{P}}
\newcommand{\tps}[1]{
\tpsSymbol{\kern-1.1pt}{(#1)}}

\newcommand{\curva}{\CurvesSymbol}
\DeclareMathOperator{\CurvesSymbol}{{\mathsf{c}}}

\DeclareMathOperator{\Fatou}{
{{\sofor}_{\mathsf{Fatou}}}
}

\newcommand{\bchar}[1]{\mathsf{#1}}

\DeclareMathOperator{\FatousetSymbol}{\mathsf{Fatou}}
\newcommand{\Fatouset}[1]{ % 1 function 
\FatousetSymbol[{#1}]
}

\DeclareMathOperator{\cvsetSymbol}{\mathtt{conv}}
\newcommand{\cvset}[2]{ % 1 function 2 set
\cvsetSymbol[#1,#2]
}

\DeclareMathOperator{\relFatouSymbol}{\FatousetSymbol}
\newcommand{\relFatou}[2]{ % 1 function 2 set
\relFatouSymbol_{#2}(#1)
}

\DeclareMathOperator{\divsetSymbol}{\mathtt{DIV}}
\newcommand{\divset}[2]{ % 1 function 2 set
\divsetSymbol[#1,#2]
}

\DeclareMathOperator{\HSymbol}{H}
\newcommand{\Hinfty}{\HSymbol^{\infty}(\udone)}
\newcommand{\Hinftyudone}{\Hinfty}

\newcommand{\hinftyd}{\hSymbol^{\infty}(\domain)}

\newcommand{\hd}{E_{n}^{+}}

\DeclareMathOperator{\hSymbol}{h}
\newcommand{\hinftyudone}{\hSymbol^{\infty}(\udone)}

\newcommand{\distanza}[3]{%1=first set, 2=second set, 3=radius
\dist\left[
A_1\cap\partial B(w,r),A_1\cap\partial B(w,r)
\right]
}

\newcommand{\eqdef}{\overset{\mathrm{def}}{=\joinrel=}}

\DeclareMathOperator{\TT}{\mathbb{T}}
\DeclareMathOperator{\NN}{{\mathbb{N}}}
\DeclareMathOperator{\ZZ}{{\mathbb{Z}}}
\DeclareMathOperator{\RR}{\mathbb{R}}
\DeclareMathOperator{\dist}{dist}

\DeclareMathOperator{\uhs}{%BoundaryOfTheUnitDisc
\mathbb{U}
}

\DeclareMathOperator{\CC}{%BoundaryOfTheUnitDisc
\mathbb{C}
}

\newtheoremstyle{slthmstyle}  % name of the style to be used
   {}       % measure of space to leave above the theorem. E.g.: 3pt
   {}       % measure of space to leave below the theorem. E.g.: 3pt
   {\slshape}   % name of font to use in the body of the theorem
   {}        % measure of space to indent
   {\bfseries}  % name of head font
   {}   % punctuation between head and body
   {2mm}       % space after theorem head
   {}           % Manually specify head
\theoremstyle{slthmstyle}
\newtheorem{theorem}{Theorem}[section]
\newtheorem*{corollary*}{Corollary.}

\newtheorem{lemma}[theorem]{Lemma}
\newtheorem{corollary}[theorem]{Corollary}
\newtheorem{proposition}[theorem]{Proposition}
\theoremstyle{definition}

\newtheorem{example}[theorem]{Example}

\newtheorem{remark}[theorem]{Remark}
\newtheorem{claim}[theorem]{Claim}
\theoremstyle{definition}

\newtheorem*{FatouThm}{Fatou's Theorem (1906)}
\newtheorem*{CorollaryL1927}{LT (1927)}
\newtheorem*{CorollaryL2006}{TLT (2006)}
\newtheorem*{Littlewood}{Littlewood's Theorem (1927)}

\newtheorem*{ATLT}{Aikawa's Theorem of Littlewood Type (1991)}
\newtheorem*{DBSSW}{Theorem LT (2006)}
\newtheorem*{DBSSWbis}{A Theorem of Littlewood Type (2006)}
\newtheorem*{Rudin'sTheorem}{Rudin's Theorem (1979)}
\newtheorem*{Rudin'sVOIFunction}{A Very Oscillatory Inner Function (1979)}
\newtheorem*{NagelSteinTheorem}{Nagel and Stein's Theorem (1984)}
\newtheorem*{NagelSteinPhenomenon}{The Nagel--Stein Phenomenon (1984)}
\newtheorem*{independenceresult}{An Independence Result (2006)}

\numberwithin{equation}{section}

\title{A New Littlewood-Type Theorem for Bounded Holomorphic Functions in the Unit Disc}
\author{Fausto Di Biase  \\
Dipartimento di Economia\\
Universit\`a ``G.\! D'Annunzio'' di Chieti-Pescara,\\
Viale Pindaro, 42, I-65127, Pescara, Italy\\
email: {\tt fausto.dibiase@unich.it}\\
\\	Haguma Gratien \\
Department of Mathematics\\ University of Rwanda\\
KN 67 Nyarugenge, Kigali, 3900 Kigali, Rwanda\\
email: {\tt hagugrat@gmail.com}\\
\\	Olof Svensson \\
Department of Science and Technology\\ Link\"oping University\\
SE-60174 Norrk\"oping,  Sweden\\
email: {\tt olof.svensson@liu.se}
}

\date{\today}
\begin{document}

\maketitle

\begin{abstract}
We prove a new theorem \textit{of Littlewood type} for bounded holomorphic functions in the unit disc $\udone$, i.e., we show that these functions 
\textit{do not all admit} almost everywhere boundary values along \textit{certain} 
systems of tangential approach regions. 
The novelty of our theorem, in comparison with the previously known results of this kind, is that it  also applies to tangential approach regions that are \textit{sequential}. 
Indeed, while in the previous results of this kind, the tangential approach regions were required to 
be \textit{curvilinear} (Littlewood 1927), (Aikawa 1990), 
or 
at least to possess a certain topological property that 
\textit{excludes} the possibility that they may be \textit{sequential} (Di Biase Stokolos Svensson Weiss 2006), 
our result is the first of this type
that also  applies  
to tangential approach regions that are \textit{sequential}. 

In order to achieve our result, we have identified a new class of 
approach regions, called \textit{projectively adjacent}, that is not defined 
in topological terms and not depend on their 
\textit{continuous} or \textit{discrete} nature. 
Indeed, this class is so broad that it contains sequential approach regions and 
curvilinear ones, as well. 

Our result can be better appreciated if we recall that 
Nagel and Stein (1984), elaborating results of 
Rudin (1979) and Nagel, Rudin and Shapiro (1982), proved 
the existence 
of translation invariant systems of tangential 
\textit{and sequential}
approach regions 
in $\hd$
along which \textit{all} bounded holomorphic functions in the unit disc 
\textit{converge a.e.} to 
their nontangential boundary values. 

\end{abstract}

{\footnotesize

\paragraph{Keywords.} Fatou type theorems, Littlewood type theorems, 
bounded holomorphic functions, approach regions, almost everywhere convergence, tangential approach regions, nontangential approach regions, very oscillatory approach regions, 
Littlewood sets, the Nagel--Stein phenomenon, projectively adjacent approach regions. 
MSC 31B25, 42B30.}

\section{Our Main Result and Its Significance}

We prove a new theorem of \textbf{Littlewood type} for bounded holomorphic functions on 
the unit disc $\udone$ in $\CC$.  
A \textit{Littlewood type theorem} for bounded holomorphic functions 
is a statement of a “negative” type: It says that these functions 
\textit{do not all admit} almost everywhere boundary values along \textit{certain} systems of \textit{tangential} approach regions. 
The following notation is used.  

The Euclidean length of $x\in\CC$ is $|x|$, and 
the (Euclidean) ball in $\CC$ of 
center $x\in\CC$ and radius $r>0$ is 
$\ball{x}{r}\eqdef\{y\in\CC:|y-x|<r\}$. Hence 
$\udone=\ball{0}{1}$. 
The boundary of $\udone$ in $\CC$ is denoted by 
$\botud\equiv\{x\in\CC:|x|=1\}$.

The Banach space of real-valued and 
harmonic functions 
defined on $\udone$ and bounded therein is denoted by 
$\hinftyudone$, and  
$\Hinftyudone$ is the Banach algebra of holomorphic functions defined on 
$\udone$ and bounded therein.

An \textbf{approach region in
$\udone$
ending at}
$\bpoint$ is a subset $\mathsf{D}\subset\udone$ 
such that $\mathsf{D}\cap\ball{\bpoint}{r}\not=\emptyset$
$\forall r>0$. 

If $\angolo\in\NN$ and $\mathsf{D}$ is an approach region 
in $\udone$ ending at $\bpoint$, we say that 
$\mathsf{D}$ is \textbf{$\angolo\,$-projectively adjacent at $\bpoint$} if 
$\forall r>0$ 
$\exists$ $\theta_1,\theta_2\in\RR$, 
$\theta_1<\theta_2$,  
and 
(a) 
$\bpoint=e^{i\theta_k}$
for either
$k=1$ or $k=2$, and 
(b)
$\forall$ $\theta\in(\theta_1,\theta_2)$
there exists 
$\dpoint\in\mathsf{D}\cap\ball{\bpoint}{r}$
such that 
$|e^{i\theta}-\dpoint|<(1+\angolo)(1-|\dpoint|)$. 
\textit{We will show momentarily the geometric meaning of this notion.} 

A \textbf{system of approach regions in $\udone$} is a 
family 
$\af\equiv{\{\afn_{\bpoint}\}}_{\bpoint\in\botud}$ where, 
$\forall \bpoint\in\botud$, 
$\afn_{\bpoint}$
is an approach region in $\udone$ ending at $\bpoint$. 
The set of all systems of 
approach regions in $\udone$ is denoted $\bafz$. 
If $\af\in\bafz$, $\af$ is \textbf{regular} if, 
for each open $U\subset \udone $, the set 
$\left\{ \bpoint\in\botud: U \cap \afn_{\bpoint} \neq \emptyset  \right\}$
 is a measurable subset of $\botud$.  

If $\af\equiv{\{\afn_{\bpoint}\}}_{\bpoint\in\botud}$ is a system of approach regions in $\udone$, 
$\bpoint\in\botud$, and $r>0$, then we set 
$\grigio{\ballrel{\bpoint}{r}\eqdef \afn_{\bpoint}\cap \ball{\bpoint}{r}}$, and 
$\af$ is called 
a \textbf{system of tangential approach regions in $\udone$} if 
$\inf_{r>0}\sup_{\dpoint\in\ballrel{\bpoint}{r}}
\frac{1-|\dpoint|}{|\dpoint-\bpoint|}=0$, $\forall\bpoint\in\botud$. 
The set of all systems of 
tangential approach regions in $\udone$ is denoted $\baftz$. 
Hence $\baftz\subset\bafz$. 

\begin{theorem}
Let $\af$ be a system of tangential approach regions in $\udone$, and assume that
\begin{itemize}
\item $\af$  \textup{\textbf{ is regular}};
\item $\forall$ $\bpoint\in\botud$, 
$\exists$ $\angolo=\angolo_{\bpoint}\in\NN$ such that 
\noindent\textup{\textbf{$\afn_{\bpoint}$ is $\angolo\,$-projectively adjacent at $\bpoint$}}. 
\end{itemize}
Then there exists a bounded holomorphic function $\dfunction$ on $\udone$ such that  
\begin{equation}
\text{for a.e.}\bpoint\in\botud, 
\lingvalue{\dfunction}{\afn_{\bpoint}}{\bpoint}
\text{ does not exist}
\label{e:lvdne} 
\end{equation}
\label{thm:maintheorem}
\end{theorem}

In order to see  how Theorem~\ref{thm:maintheorem} differs from 
the previously known results of this kind, the following notions and facts will be helpful. 

$\bafpaz$ is the set of all 
systems of \textbf{projectively adjacent} approach regions in $\udone$, i.e., 
the set of all $\af\in\bafz$ such that $\forall$ 
$\bpoint\in\botud$ 
$\exists$ 
$\angolo\in\NN$ such that 
$\afn_{\bpoint}$ is $\angolo\,$-projectively adjacent at $\bpoint$.

$\bafcz$ is the set of all 
systems $\af$ of \textbf{curvilinear} approach regions in $\udone$, i.e., such that
$\forall$ $\bpoint\in\botud$ 
$\exists$ $r_{\bpoint}>0$ and 
there exists a continuous function 
$\varphi_{\bpoint}:[0,1)\to\udone$ 
with 
$\displaystyle{\lim_{s\uparrow1}\varphi_{\bpoint}(s)=\bpoint}$ and 
 such that 
$\displaystyle{\ballrel{\bpoint}{r_{\bpoint}}=\ball{\bpoint}{r_{\bpoint}}
 \cap\{\varphi_{\bpoint}(s):0\leq s<1\}\}}$. 

$\bafsz$ is the set of all 
systems $\af$ of \textbf{sequential} approach regions in $\udone$, i.e., such that 
$\forall$ 
$\bpoint\in\botud$ 
$\exists$ 
$r_{\bpoint}>0$ 
and there exists a sequence 
$\varphi_{\bpoint}:\NN\to\udone$
such that 
$\displaystyle{\lim_{j\to+\infty}\varphi_{\bpoint}(j)=\bpoint}$
and
$\displaystyle{\ballrel{\bpoint}{r_{\bpoint}}=\ball{\bpoint}{r_{\bpoint}}
 \cap\{\varphi_{\bpoint}(j):j\in\NN\}}$.

Then 
\begin{equation}
\text{(i)}\quad\bafpaz\cap\bafcz\not=\emptyset;\quad \text{(ii)}
\quad\bafpaz\cap\bafsz\not=\emptyset;
\label{e:good}
\end{equation}
and
$$
\bafsz\cap\baftz\cap\bafpaz\not=\emptyset\text{ and }
(\bafsz\cap\baftz)\setminus\bafpaz\not=\emptyset
$$
In particular, on the one hand, in the present work we introduce 
a new class of systems of approach regions 
(the set 
$\bafpaz\subset\bafz$ of all systems of projectively adjacent approach regions), 
which is 
does not depend on the \textit{continuous} or \textit{discrete} nature of the 
approach regions involved, and indeed 
is \textit{transversal}, so to speak, 
to the properties of being \textit{sequential} or  \textit{curvilinear}, as shown in~\eqref{e:good}. 
On the other hand, while the previously known results of Littlewood type 
  apply to systems of tangential approach regions that are \textit{curvilinear} \cite{Littlewood1927} 
and 
\cite{Zygmund1949} 
(and \cite{Aikawa1991}, for Euclidean half-spaces), 
or share with them a topological property that \textit{excludes} the possibility 
that they may be sequential 
\cite{DiBiaseStokolosSvenssonWeiss2006}, 
our result \textit{also} applies to tangential approach regions that are sequential. 

Our result can be better appreciated if we recall that 
A. Nagel and E.M. Stein \cite{Nagel--Stein1984}, elaborating a result of 
W. Rudin \cite{Rudin1979bis} and A. Nagel, W. Rudin and J.H. Shapiro \cite{NagelRudinShapiro}, proved 
the existence 
of rotationally invariant systems of tangential 
\textit{and sequential}
approach regions 
in $\udone$ (and in higher-dimensional Euclidean half-spaces as well)
along which \textit{all} bounded harmonic functions \textit{converge a.e.} to 
their nontangential boundary values. 

The following notions will be useful to  illustrate the geometric meaning of the notion of projective adjacency.

If $X$ is a set then $\tps{X}$ is the collection of all subsets of $X$. 

If $\af\in\bafz$ 
and $U\subset\udone$ then the set 
$\{\bpoint\in\botud:\afn_{\bpoint}\cap{}U\not=\emptyset\}$
is denoted by 
$\af^{\ast}(U)$ and is called the \textbf{$\af$-shadow of} $U$. 
The map  $\af^{\ast}:\tps{\udone}\to\tps{\botud}$ thus defined is called 
the $\af$-\textbf{projection}.
It follows that $\af\in\bafz$ 
is regular if and only if 
\textit{the projection associated to $\af$ maps each open subset of $\udone$ to a measurable 
subset of $\botud$}. 

If $\angolo\in\NN$ and $\bpointv\in\botud$ 
then 
$\Gamma_{\angolo}(\bpointv)\eqdef\{\dpoint\in\udone: 
(1-|\dpoint|)/(|\bpointv-\dpoint|)>1/(1+\angolo)\}$. 
%\begin{equation}
%\Gamma_{\angolo}(\bpointv)\eqdef\{\dpoint\in\udone: 
%(1-|\dpoint|)/(|\bpointv-\dpoint|)>1/(1+\angolo)\};
%\label{e:standard}
%\end{equation}
$\Gamma_{\angolo}(\bpointv)$ is called the 
\textbf{nontangential approach region at} $\bpointv$ of \textit{width} $\angolo$; 
cf. \cite{SteinWeiss1971}. Observe that 
${\boldsymbol{\Gamma}}_{\boldsymbol{\angolo}}={\{\Gamma_{\angolo}(\bpointv)\}}_{\bpointv\in\botud}$ 
is a system of approach regions in $\udone$.  

If $\dfunction:\udone\to\CC$, the \textbf{Fatou set of} $\dfunction$ is defined as 

$$
\Fatouset{\dfunction}
\eqdef
\{\bpoint\in\botud: 
\exists\ell\in\CC
\text{ such that }
\lim_{\stackrel{\text{nt}}{z\to\bpoint}}\dfunction(z)=\ell\},
$$ 
where the notation 
$\grigio{\displaystyle{\lim_{\stackrel{\text{nt}}{z\to\bpoint}}\dfunction(z)=\ell}}$ means 
that $\forall\angolo\in\NN$ and $\forall\epsilon>0$ 
$\exists\delta>0$ such that if $z\in\ball{\bpoint}{\delta}\cap\Gamma_{\angolo}(\bpoint)$ 
then $|\dfunction(z)-\ell|<\epsilon$. 
If $\dfunction\in\hinftyudone$ then $\Fatouset{\dfunction}$
has 
full measure in $\botud$ (see \cite{Zygmund1988} and references therein), 
and the  
\textbf{nontangential boundary function of} $\dfunction$ is 
the function 
$\ntbv{\dfunction}:\Fatouset{\dfunction}\to\RR$, defined by 
$$
\ntbv{\dfunction}(\bpoint)\eqdef
\lim_{\stackrel{\text{nt}}{z\to\bpoint}}\dfunction(z). 
$$
If $\bpoint\in\botud$ and $0<r<\pi$ the set 
\begin{equation}
\grigio{\acal{\bpoint}{r}\eqdef\{\bpoint{}e^{i\theta}:|\theta|<r\}
}
\label{e:intcentered} 
\end{equation}
is called 
an \textbf{arc in $\botud$} of \textbf{center $\bpoint$} and \textbf{endpoints} 
$\bpoint{}e^{ir}$ and $\bpoint{}e^{-ir}$, and 
$\bpoint{}e^{ir}$ 
($\bpoint{}e^{-ir}$ )
is 
the \textbf{right-endpoint} 
(the \textbf{left-endpoint} )
of $\acal{\bpoint}{r}$. 
The collection of all arcs in $\botud$ is denoted $\soba$. 
We also set 
$\grigio{\racal{\bpoint}{\theta}\equiv\{\bpoint{}e^{i\beta}:0<\beta<\theta\}}$. 

If ${\mathsf{D}}$ is an approach regions in $\udone$ 
ending at $\bpoint$, then  ${\mathsf{D}}$ is 
$\angolo\,$-projectively adjacent at $\bpoint$ if and only if 
\textit{$\forall$
$r>0$ 
the $\boldsymbol{\Gamma_b}$-shadow of ${\mathsf{D}}\cap\ball{\bpoint}{r}$ contains 
an arc of which $\bpoint$ is an endpoint}.

The following background results also help us assess the novelty of this work. 
 
In 1949, A. Zygmund gave a real-variable proof of Littlewood's theorem, that may be extended to higher-dimensional Euclidean half-spaces $\hd$
(see \cite{SteinWeiss1971} for background)
 so as to yield the following result: If 
$\af$ is translation invariant system of approach regions and, for each $\bpoint\in\bdryhd$, 
$\af(\bpoint)$ is a hypersurface of dimension $n-1$ tangential to $\bdryhd$, then 
there exists a bounded harmonic function $\dfunction$ on $\hd$ 
such that~\eqref{e:lvdne} holds.

In 1991, H. Aikawa proved that 
a deep result, a theorem of Littlewood type for $\hd$ where the 
tangential approach regions are supposed to be curvilinear (see \cite{Aikawa1991} 
for a precise statement). 

In 2006, the statement 
\flqq\textit{there is no $\af\in\baftz\cap\bafcz$ such that each 
$\dfunction\in\Hinftyudone$
converges to its nontangential boundary values almost everywhere along $\af$}\frqq\,
was proved to be \textit{independent of} ZFC \cite{DiBiaseStokolosSvenssonWeiss2006}. 
In other words, it is neither possible to prove it, nor to disprove it.

The property of projective adjacency, 
introduced in the present work, should be compared 
to the notion 
of being \textit{attached}, 
introduced in \cite[p. 48]{DiBiaseStokolosSvenssonWeiss2006} 
and therein denoted (c$\star$). If we denote by $\bafaz$ the collection of all systems of 
\textit{attached} approach regions in $\udone$, then  
\begin{equation}
\text{(i)}\quad\bafcz\subsetneq\bafaz\quad\text{ but }\quad \text{(ii)}\quad\bafaz\cap\bafsz=\emptyset 
\label{e:attachednew} 
\end{equation}
It is useful to compare \eqref{e:attachednew} (ii) with~\eqref{e:good} (ii).

\subsection{Additional notation and background}

The function 
$\tau : \botud \times \udone \to    (0, 1]$ 
defined by 
$\tau(\bpoint,\dpoint) \eqdef (1-|\dpoint|)/(|\bpoint-\dpoint|)$, 
for 
$\bpoint\in\botud$ and $\dpoint\in\udone$, 
is the 
\textbf{normalized distance to the boundary}. 

A \textbf{null set} in $\botud$ is a measurable subset 
$\Ub\subset\botud$ with $|\Ub|=0$,
 and we set $\nullsets\eqdef\{\Ub\in\tps{\botud}:
\Ub\text{ is a null set in } \botud\}$. 
If $\Sb$ and $\Qb$ are subsets of $\botud$ we say that 
$\Sb$ is \textbf{a.e. contained in} $\Qb$, and write
$\grigio{\aesubset{\Sb}{\Qb}}$
if $\Sb\setminus\Qb\in\nullsets$: 
This means that almost all of $\Sb$ is a subset of $\Qb$. Observe that 
the sets $\Sb$ and $\Qb$ need not be measurable
and that 
\begin{equation}
\aesubset{\Ab\setminus\Bb}{\Cb}
\iff
\aesubset{\Ab\setminus\Cb}{\Bb}.
\label{e:aesubset} 
\end{equation}
We say that $\Sb$ and $\Qb$ are \textbf{almost everywhere equal}, 
and write 
$\grigio{\aeequal{\Sb}{\Qb}}$
if $\aesubset{\Sb}{\Qb}$ and $\aesubset{\Qb}{\Sb}$.
A set $\Sb\subset\botud$ 
has 
{\bf full measure} 
if $\aeequal{\Sb}{\botud}$. 
A property is said to hold 
\textbf{a.e.}  
if the set of points 
in $\botud$ 
for which it holds has full measure. 
A set $\Sb\subset\Qb$ 
has 
\textbf{full measure in} $\Qb$ 
if $\aeequal{\Sb}{\Qb}$. 

Recall that $z\in\CC$ has \textbf{polar coordinates} $(r,\theta)$ if 
$z=r\exp(i\theta)$, where $r\geq0$, $\theta\in\RR$. 
Points in a neighborhood of $\botud$
may be 
 conveniently parameterized by 
a natural variant of the polar coordinates: We say that $z$ has 
\textbf{boundary coordinates} $(\theta,\delta)$ if 
$z=(1-\delta)e^{i\theta}$. Hence 
$\delta$ is the \textit{signed distance to the boundary}, since $\delta<0$ if 
$z$ lies outside of $\udone$.  
For example, if $O$ is the neighborhood of $\bpoint=1$ given by 
\begin{equation}
O\eqdef\{(1-\delta)e^{i\theta}:|\delta|<\frac{1}{4},|\theta|<\frac{1}{4} \}
\label{e:openset} 
\end{equation}
and we denote by $S$ the square ${(-\frac{1}{4},\frac{1}{4})}^2\eqdef
\{(\theta,\delta):|\theta|<\frac{1}{4},|\delta|<\frac{1}{4} \}$, then 
for the function  
$p:{S}\to{O}$ defined by 
$p(\theta,\delta)\eqdef(1-\delta)e^{i\theta}$ we have 
 the following approximation, which is rough but useful. 
If  $(\theta_k,\delta_k)\in S$, $k=1,2$,  then 
\begin{multline}
\frac{3}{64}\left[{(\theta_1-\theta_2)}^2+{(\delta_1-\delta_2)}^2\right]
\leq 
{|p(\theta_1,\delta_1)-p(\theta_2,\delta_2)|}^2
\\
\leq \frac{125}{64}\left[
{(\theta_1-\theta_2)}^2+{(\delta_1-\delta_2)}^2
\right]
\label{e:estimate}
\end{multline}
The proof of~\eqref{e:estimate} is achieved by 
a laborious but essentially straightforward calculation based on 
Taylor's second order approximation of $\cos x$.

\paragraph{Carleson tents and the Zygmund map.}

The \textbf{Carleson tent above} an arc 
$\acal{\yb}{\theta}\subset\botud$ is defined by 
\begin{equation}
\Delta(\acal{\yb}{\theta})\eqdef
\udone\cap\ball{\yb}{|\yb-\yb{}e^{i\frac{\theta}{2}}|}.
\label{e:Delta} 
\end{equation}
The \textbf{Zygmund map} is the map 
$\Zb:\tps{\botud}\to\tps{\botud}$
defined as follows: If $\Vb\subset\botud$ then $\Zb(\Vb)\subset\botud$ 
is the set of points of $\botud$ with the following properties:
\begin{equation}
\textup{(i)}\,
\bpoint\in\botud\setminus\Vb;\quad
\textup{(ii)}\,
\forall\epsilon>0\,
\exists 
\text{ arc } \Jb_{\epsilon}\subset\Vb\cap\ball{\bpoint}{\epsilon}
\text{ with }
\af(\bpoint)\cap\Delta(\Jb_{\epsilon})\not=\emptyset
\label{e:Zygmund} 
\end{equation}

\paragraph{Poisson integrals and the basic Carleson tent estimate.}

If $\bfunction:{\botud}\to{\RR}$ is Lebesgue integrable then 
the \textbf{Poisson integral} of 
$\bfunction$ is denoted by $\PoissonI{\bfunction}$. 
If $\Sb\subset\botud$ then 
$\cfunction{\Sb}:{\botud}\to{\RR}$ 
is defined by 
$$
\cfunction{\Sb}(\bpoint)=1\,
\text{ if }
\bpoint\in\Sb,
\,\,
\cfunction{\Sb}(\bpoint)=0
\text{ if }
\bpoint\not\in\Sb.
$$
The following result is well-known (cf.\ \cite{Calderon1950}, 
\cite{SteinWeiss1971}). 
The constant $c_0>0$ that appears therein is called 
the \textbf{Carleson tent constant}. 
\begin{lemma} 
There exists a constant $c_0\in(0,+\infty)$ with the following property: 
For each arc $\Jb\subset\botud$ and each $\dpoint\in\Delta(\Jb)$, 
$\PoissonI{
\cfunction{\Jb}}
(\dpoint)\geq{}c_0$.  
\label{l:CarlesonTent}
\end{lemma}
\begin{lemma}
If $\Vb \subset \botud$ is open then, for all $\bpoint \in \Zb(\Vb)$,
$$
\limsup_{\stackrel{\dpoint\to\bpoint}{\dpoint\in\af(\bpoint)}}  
\PoissonI{\cfunction{\Vb}}(\dpoint)\geq c_0 
$$
\end{lemma}
\begin{proof}
It follows from Lemma~\ref{l:CarlesonTent}, since if $\Jb \subset \Vb$ then 
$\PoissonI{\cfunction{\Jb}}\leq\PoissonI{\cfunction{\Vb}}$.
\end{proof}

\section{Proof  of Theorem~\ref{thm:maintheorem}}
\label{s:proofofthemaintheorem}

If $\dfunctionharm \in \hinftyudone$ and $\bpoint \in \botud$,  define
$\displaystyle{\osc{\dfunctionharm}{w}\eqdef 
\limsup_{\stackrel{\dpoint\to\bpoint}{\dpoint\in\af(\bpoint)}} 
\dfunctionharm(\dpoint)
-
\liminf_{\stackrel{\dpoint\to\bpoint}{\dpoint\in\af(\bpoint)}} 
\dfunctionharm(\dpoint)}$, and observe that, 
in order to prove Theorem~\ref{thm:maintheorem}, it suffices to prove 
the following statement
\begin{equation}
\flqq\exists\,
\dfunctionharm\in\hinftyudone
\text{ with } 
\dfunctionharm>0
\text{ and } 
\osc{\dfunctionharm}{\bpoint}>0
\text{ for a.e. } 
\bpoint\in\botud
\frqq
\end{equation}
Indeed, 
if $\dfunctionharmc$ is the harmonic conjugate to $\dfunctionharm$ then 
$\dfunctionhol\eqdef{}e^{-\dfunctionharm-i\dfunctionharmc}$ has the required properties. 
Since $\af$ satisfies the hypotheses of Theorem~\ref{thm:maintheorem}, 
%we may assume, without loss of generality, that 
%\begin{equation}
%\angolo_{\bpoint}\geq 10\quad\forall\bpoint\in\botud
%\label{e:largeangle} 
%\end{equation}
consider the following sequence of everywhere defined functions 
$\bfunction_j: \botud \to (0, \infty)$ gauging
the order of tangency of $\afn_{\bpoint}$ at various points: 
\begin{equation}
\bfunction_j(\bpoint) \eqdef \sup\left\{
\tau(\bpoint, \dpoint) :
\dpoint\in\ballrel{\bpoint}{(22\pi)/(j10)}
\right\} 
\label{e:fn}
\end{equation}
Consider the sequence $v:{\NN}\to{\NN}$ whose values are, in the 
natural order:
$$
2,3,2,3,4,2,
3,4,5,2,3,4, \ldots
$$
so that for each $\bsegnato\in\NN\setminus\{1\}$ there are infinitely many values of $j\in\NN$ such that 
$v_j=\bsegnato$. 
Since $\afn_{\bpoint}$ ends tangentially at $\bpoint$, 
for each $\bpoint\in\botud$ the sequence 
${\{\bfunction_n(\bpoint)\}}_{n\in\NN}$ decreases to $0$. Moreover, since 
$\af$ is regular, for 
each $n\in\NN$ the function $\bfunction_n:\botud\to(0,1]$ is measurable. 
Egorov's Theorem \cite{Royden1968bis} implies that  
for each $j \in \mathbb{N}$  then there is a set $\Cb_j \subset \botud$ whose Lebesgue measure is greater than $2\pi-2^{-j}$ and such that the sequence $\left\{f_n\right\}$ converges uniformly to 0 on $\Cb_j$. 
Observe that 
if the sequence $\{\bfunction_n\}$ converges uniformly to $0$ on a set $\Ab$ and on a set $\Bb$, then it
converges uniformly to $0$ on 
 $\Ab\cup \Bb$. Hence 
we may assume that 
\begin{equation}
\Cb_j \subset \Cb_{j+1}
\quad\forall\, j\in\NN
\label{e:Cisincreasing} 
\end{equation}
It follows that there exists a sequence $\phi:{\NN}\to{\NN}$ 
such that, for each $j\in\NN$ 
\begin{equation}
\textup{(i)}
\,\,
\sup \left\{f_{{\phi_j}}(\vb): \vb \in \Cb_j\right\} 
< c\frac{10}{22}\frac{2^{-j}}{v_j},\,\,
\textup{(ii)}
\,\,
\phi_j>j,
\,\,
\textup{(iii)}\,
\phi_{j+1}>\phi_j
\label{e:aenne} 
\end{equation}
where $c$ is a constant that will be chosen momentarily. 
Define 
\begin{equation}
\Cb\eqdef\bigcup_{j\in\NN}\Cb_j
\label{e:dofC} 
\end{equation}
Observe that $\aeequal{\Cb}{\botud}$. 

\paragraph{Basic constructions}

The following devices will be employed in the proof of 
Theorem~\ref{thm:maintheorem}. 
Let $\bpoint\in\botud$ and $|\theta|\leq 2\pi$. 
Recall that $\acal{\bpoint}{\theta}$, defined in~\eqref{e:intcentered}, 
is the arc in $\botud$ 
centered at $\bpoint$ of arc-length $\theta$.
If $\Yb\subset\botud$ 
we define
$$
\Acal{\Yb}{\theta}\eqdef\bigcup_{\bpoint\in\Yb}\acal{\bpoint}{\theta}
$$
If $n\in\NN$ and $n\geq1$, we define 
\begin{equation}
\Ub_n\eqdef\{e^{i2\pi{}p/n}:p=0,1,\ldots,n-1\}.
\label{e:definitionofy}
\end{equation}
Then 
$|\Acal{\Ub_n}{2\pi/n}|=2\pi$, 
since $\Acal{\Ub_n}{2\pi/n}$ is the union of $n$ disjoint arcs
whose union is equal to $\botud\setminus\{e^{i\pi/n}e^{ip2\pi/n}:p=0,1,\ldots{}n-1\}$. 
Indeed,  
$$
\Acal{\Ub_n}{2\pi/n}=
\bigcup_{\bpoint\in\Ub_n}\acal{\bpoint}{2\pi/n}=
\bigcup_{p=0}^{n-1}\acal{e^{i2\pi{}p/n}}{2\pi/n}
$$
and the arcs appearing in the union are disjoint. 
Define 
\begin{equation}
\Ob_n  \eqdef \Acal{\Ub_{{\phi_n}}}{2^{-n}2\pi/{\phi_n}}
\label{e:definitionofOn}
\end{equation}
Observe that 
$\left|\Ob_n\right|=2^{-n}2\pi$
and hence 
$\lim_{j\to+\infty}\left|\bigcup_{k\geq{}j}\Ob_k\right|=0$. 
Define 
\begin{equation}
\Vb_j\eqdef\bigcup_{k\geq{}j}\Ob_k
\label{e:Vn} 
\end{equation}
The set $\Vb_j$ is open and dense in $\botud$, and, if 
$j$ is large, it has small measure. 

%\subsection{Basic lemmata}
%\label{s:basiclemmata}

\paragraph{Basic lemmata.}
\begin{lemma} If $\af$ satisfies the hypotheses of 
\textup{Theorem~\ref{thm:maintheorem}} then 
each $n\in \mathbb{ N}$,
$\Cb \setminus \Vb_n \subset \Zb(\Vb_n)$. 
\label{l:bl1}
\end{lemma}
\begin{proof}
Before we delve into the proof, we would like to explain why it is more 
involved than one would perhaps think. 
It is indeed true that, given $\bpoint,\bpointv\in\botud$, with 
$\bpoint\not=\bpointv$, and given $b\in\NN$ and $r>0$ there exists
$n=n(\bpoint,\bpointv,b,r)\in\NN$ such that 
$\Gamma_{b}(\bpointv)\setminus\Gamma_{n}(\bpoint)\subset\ball{\bpointv}{r}$: 
However, in this statement, $n$ depends upon $\bpoint$, 
$\bpointv$, 
$b$,
and $r$, 
while, in the statement that we need to prove, the variables appear in a different order 
and hence the dependence of $n$ on the other variables is not admissible. 
Let $\bopoint\in\Cb\setminus\Vb_n$. Then there exists $\dsegnato\in\NN$ with 
$\bpoint\in\Cb_{\dsegnato}$, and it follows that 
\begin{equation}
\grigio{j\geq \dsegnato\implies \bopoint\in\Cb_j
\text{ and hence }
\eqref{e:aenne}
\text{ implies that }
f_{\phi_j}(\bopoint)\leq{}c\frac{10}{22}\frac{2^{-j}}{v_j}
}
\label{e:fromd}
\end{equation}
Moreover, 
\begin{equation}
\forall{}j\geq n, \,\,\bopoint\not\in \Ob_j
\label{e:fromdBIS}
\end{equation}
and there exists $\bsegnato\in\NN$ such that 
\begin{equation}
\afn_{\bopoint}
\text{ is }
\bsegnato\text{-projectively adjacent},\,\,\bsegnato\geq10
\label{e:bsegnato}
\end{equation}
Let $\epsilon>0$. 
Since 
$\afn_{\bpoint}$ ends tangentially at $\bpoint$, 
there exists $r>0$ such that 
\begin{equation}
\sup\{\tau(\bpoint,\dpoint):\dpoint\in
\ballrel{\bpoint}{r}
\}<\frac{1}{10}\frac{1}{{(1+\bsegnato)}^2}
\text{ and }
r<\epsilon\cdot{}2^{-10}
\label{e:epsilon'} 
\end{equation}
Then~\eqref{e:bsegnato} implies 
that  
\begin{equation}
\Gamma_{\bsegnato}^{\ast}(\ballrel{\bpoint}{r})
\text{ contains an arc $\Jb$ of which }
\bpoint
\text{ is an endpoint }
\label{e:padj} 
\end{equation}
%\begin{equation}
%\Gamma_{\bsegnato}^{\ast}(\ballrel{\bpoint}{r})\leadsto\bpoint 
%\label{e:padj} 
%\end{equation}
We may assume, without loss of generality, 
that 
$\bpoint$ is a left-endpoint of 
$\Jb$, i.e., that 
there exists $\theta>0$ such that 
\begin{equation}
\textup{(i) }
\,
\racal{\bpoint}{\theta}
\equiv
\{\bpoint{}e^{i\beta}:0<\beta<\theta\}
\subset\Gamma_{\bsegnato}^{\ast}(\ballrel{\bpoint}{r}),
\,
\textup{ (ii) }
\,
\theta<r\cdot{}2^{-10}.
\label{e:theta} 
\end{equation}
Select $j\in\NN$ such that 
\begin{equation}
\textup{(i)}\,\, j>n,
\quad
\textup{(ii)}\,\,
2\pi/{\phi_j}<\theta\cdot{2^{-10}},
\quad
\textup{(iii)}
\,\,
j>\dsegnato,
\quad
\text{ and }
\textup{(iv)}
\,\,
v_j=1+\bsegnato
\label{e:selectenne} 
\end{equation}
It follows that there exists 
$p\in\ZZ$ such that if 
$\vb\eqdef{}e^{i2\pi{}p/{{\phi}_{j}}}$ then 
\begin{equation}
\bpoint\in
\racal{\vb}{2\pi/{\phi_j}}
\equiv
\{\vb{}e^{i\alpha}:0<\alpha<{2\pi}/{\phi_j}\}
\label{e:contained}
\end{equation}
Let
$\yb\eqdef{}
e^{i2\pi{}(p+1)/{{\phi}_{j}}}
\equiv{}
\vb{}
e^{i{2\pi}/{\phi_n}}$. 
Since 
$\bpoint\not\in \Ob_j$ it follows that 
\begin{equation}
\bpoint \not\in\acal{\vb}{2^{-j}\frac{2\pi}{\phi_j}}
\text{ and }
\bpoint \not\in\acal{\yb}{2^{-j}\frac{2\pi}{\phi_j}}
\end{equation}
We claim that $\yb\in\racal{\bpoint}{\theta}$. 
Indeed, \eqref{e:contained} 
and~\eqref{e:selectenne} (ii)
imply that 
$\bpoint=\vb{}e^{i\alpha}$ with $0<\alpha<2\pi/\phi_j<\theta\cdot{}2^{-10}$
and
$\yb=\vb{}e^{i{2\pi}/{\phi_j}}=
\vb{}e^{i\alpha}e^{i({2\pi}/{\phi_j}-\alpha)}=\bpoint{}e^{i({2\pi}/{\phi_j}-\alpha)}$
with $0<{2\pi}/{\phi_j}-\alpha<{2\pi}/{\phi_j}<\theta\cdot{2^{-10}}<\theta$.
Hence 
\begin{equation}
\yb=\bpoint{}e^{i\beta}
\text{ with }
0<\beta<{2\pi}/{\phi_j}
\label{e:estimate1} 
\end{equation}
where $\beta\eqdef{2\pi}/{\phi_j}-\alpha$.
Since $0<\beta<{2\pi}/{\phi_j}<\theta/20<\theta$, 
the claim is proved. 

Since $\yb\in\racal{\bpoint}{\theta}$, \eqref{e:theta} implies that 
$\yb\in \Gamma_{\bsegnato}^{\ast}(\ballrel{\bpoint}{r})$, hence 
\begin{equation}
\text{there exists }
\zsegnato\in\ballrel{\bpoint}{r}
\text{ such that }
\zsegnato\in\Gamma_{\bsegnato}(\yb)
\label{e:projection} 
\end{equation}
Since 
$\zsegnato\in\ballrel{\bpoint}{r}$, \eqref{e:epsilon'} implies that 
$\tau(\bpoint,\zsegnato)<\frac{1}{10}\frac{1}{{(1+\bsegnato)}^2}$, i.e.
\begin{equation}
1-|\zsegnato|<|\bpoint-\zsegnato|\frac{1}{10}\frac{1}{{(1+\bsegnato)}^2}
\label{e:alzata} 
\end{equation}
Since $\zsegnato\in\Gamma_{\bsegnato}(\yb)$ it follows that 
$\tau(\yb,\zsegnato)>1/(1+\bsegnato)$, i.e. 
\begin{equation}
|\yb-\zsegnato|<(1-|\zsegnato|) (1+\bsegnato)
\label{e:alzata1}
\end{equation}
Then~\eqref{e:alzata} and~\eqref{e:alzata1} imply that 
\begin{equation}
|\yb-\zsegnato|<|\bpoint-\zsegnato| \frac{1}{10}\frac{1}{1+\bsegnato}
\label{e:alzata2}
\end{equation}
Then $|\bpoint-\zsegnato|\leq|\bpoint-\yb|+|\yb-\zsegnato|<|\bpoint-\yb|+|\bpoint-\zsegnato| 
\frac{1}{10}\frac{1}{1+\bsegnato}$, hence
\begin{equation}
|\bpoint-\zsegnato|\left[1-\frac{1}{10}\frac{1}{1+\bsegnato}\right]<|\bpoint-\yb|
\label{e:alzata3}
\end{equation}
Observe that~\eqref{e:estimate1} implies that 
\begin{equation}
|\bpoint-\yb|<2\pi/\phi_j
\label{e:alzata4}
\end{equation}
and then~\eqref{e:alzata3} implies that  
\begin{equation}
|\bpoint-\zsegnato|<\frac{1+\frac{1}{\bsegnato}}{1+\frac{9}{10\bsegnato}}
\frac{2\pi}{\phi_j}
<\frac{11}{10}\frac{2\pi}{\phi_j}
\label{e:alzata5}
\end{equation}
Hence 
\begin{equation}
\zsegnato\in\ballrel{\bpoint}{\frac{11}{10}\frac{2\pi}{\phi_j}}
\label{e:alzata7bis}
\end{equation}
and thus~\eqref{e:fn} implies that 
\begin{equation}
\tau(\bpoint,\zsegnato)<f_{\phi_j}(\bpoint)
\label{e:alzata7}
\end{equation} 
Since 
$\bpoint\in\Cb_j$, \eqref{e:aenne} 
and~\eqref{e:selectenne}
imply that 
\begin{equation}
f_{\phi_j}(\bpoint)<c\frac{10}{22}\frac{2^{-j}}{v_j}
=
c\frac{10}{22}\frac{2^{-j}}{1+\bsegnato}
\label{e:alzata8}
\end{equation} 
\eqref{e:alzata7} and~\eqref{e:alzata8} imply that 
$\tau(\bpoint,\zsegnato)=\frac{1-|\zsegnato|}{|\bpoint-\zsegnato|}
<c\frac{10}{22}\frac{2^{-j}}{1+\bsegnato}$, i.e., 
\begin{equation}
1-|\zsegnato|<c
|\bpoint-\zsegnato|
\frac{10}{22}\frac{2^{-j}}{1+\bsegnato}
\label{e:schiacciata}
\end{equation} 
and now~\eqref{e:alzata1}, 
\eqref{e:schiacciata}, and \eqref{e:alzata5}
imply that, thanks to~\eqref{e:estimate}

\begin{equation}
|\yb-\zsegnato|<c\frac{1}{2}2^{-j}\frac{2\pi}{\phi_j}
\leq{}c\sqrt{\frac{64}{3}}|\yb-\yb{}e^{i\frac{2^{-j}{2\pi}/{\phi_j}}{2}}|
=|\yb-\yb{}e^{i\frac{2^{-j}{2\pi}/{\phi_j}}{2}}|
\label{e:schiacciatabis}
\end{equation} 
where $c=\sqrt{\frac{3}{64}}$
and hence 
\begin{equation}
\zsegnato\in \Delta(\acal{\yb}{2^{-j}{2\pi}/{\phi_j}})
\label{e:schiacciatatent}
\end{equation} 
This fact concludes the proof.  
\end{proof}

\begin{lemma}
If $\af$ satisfies the hypotheses of \textup{Theorem~\ref{thm:maintheorem}}, 
$\Vb \subset \botud$ is open, and 
$|\botud \setminus \Vb| > 0$, then 
\begin{equation}
\aesubset{\botud \setminus \Vb}{
\left\{\bpoint\in\botud:\liminf_{\stackrel{\dpoint\to\bpoint}{\dpoint\in\afn_{\bpoint}}} 
(\PoissonI{\cfunction{\Vb}})(\dpoint) =0\right\}}
\label{e:C1950}
\end{equation}
\label{l:C1950}
\end{lemma}
\begin{proof}
If we denote by 
$\Qb$ the set in the right-hand side of~\eqref{e:C1950}, then we have 
to prove that $\aesubset{\botud\setminus\Vb}{\Qb}$, i.e., 
that $(\botud\setminus\Vb)\setminus\Qb\in\nullsets$. 
The proof of Lemma~\ref{l:C1950} rests on the following claim. 

\begin{claim}
For each $\delta>0$ 
there exists  $\Eb\in\tps{\botud}$ with 
$|\Eb|<\delta$ and 
\begin{equation}
\aesubset{(\botud\setminus\Vb)\setminus\Eb}{\Qb}
\label{e:claim1} 
\end{equation}
\label{claim:one} 
\end{claim}

Claim~\ref{claim:one} and~\eqref{e:aesubset} imply that 
$\aesubset{(\botud\setminus\Vb)\setminus\Qb}{\Eb}$, and hence 
$|(\botud\setminus\Vb)\setminus\Qb|<\delta$, 
and since $\delta>0$ is arbitrary, it follows that 
$(\botud\setminus\Vb)\setminus\Qb\in\nullsets$, thus completing the proof 
of~Lemma~\ref{l:C1950}.

\paragraph{Proof of Claim~\ref{claim:one}.} 
Let $\delta>0$. 
A well-known extension of Fatou's Theorem to Poisson integrals 
(see e.g. \cite{SteinWeiss1971}) says that 
$$
\aeequal{\Fatouset{\PoissonI{
\cfunction{\Vb}
}}}{\botud}
\,\,
\text{ and }
\,\,
\ntbv{(\PoissonI{\cfunction{\Vb}})}(\bpoint)=1_{\Vb}(\bpoint) 
\,\,\,\forall\bpoint\in\Fatouset{P(1_{\Vb})}.
$$ 
It follows that 
$\ntbv{(\PoissonI{
\cfunction{\Vb}})}(\bpoint)=0
\quad
\text{ for a.e. }
\bpoint\in\botud\setminus\Vb$.
We denote by $\Xb$ the subset of $\botud\setminus\Vb$ of full measure 
in $\botud\setminus\Vb$ with the property that 
$\ntbv{(\PoissonI{\cfunction{\Vb}})}(\bpoint)=0$
for all $\bpoint\in\Xb$, and define, 
for each $(j,b)\in\NN\times\NN$, the function 
$$
g_{(j,b)}:\Xb\to(0,1]
\quad
\text{ by }
g_{(j,b)}(\bpoint)\eqdef\sup\{
\PoissonI{\cfunction{\Vb}}(\dpoint):\dpoint\in\Gamma_{b}(\bpoint)\cap\ball{\bpoint}{1/j}\}.
$$
Observe that 
$\forall b,j\in\NN,\,\,
g_{(j,b)}\geq g_{(j+1,b)}
\text{ and }
g_{(j,b)}\leq g_{(j,b+1)}$ 
and that 
$$
\forall\bpoint\in\Xb,\,\lim_{j\to+\infty}g_{(j,b)}(\bpoint)=0
$$
Egorov's theorem implies that for each $k\geq1$  there 
exists a set $\Eb(k)\subset\Xb$ with 
\begin{description}
\item[(a)] $|\Eb(k)|<\delta\cdot{}2^{-k}$;
\item[(b)] the sequence 
${\{g_{(j,k)}\}}_{j\in\NN}$ converges to zero uniformly on $\Xb\setminus\Eb(k)$;
\item[(c)] the set $\Xb\setminus\Eb(k)$ is perfect.
\end{description}
Let 
$\Eb\eqdef\cup_{k\geq1}\Eb(k)$. Then 
$|\Eb|<\delta$, 
$\Xb\setminus\Eb=\cap_{k\geq1}\Xb\setminus\Eb(k)\subset\Xb\setminus\Eb(b)$ for each $b\geq1$. 
The proof of Claim~\ref{claim:one} rests on the following result. 
\begin{claim}
\begin{equation}
\aesubset{\Xb\setminus\Eb}{\Qb}
\label{e:claim2} 
\end{equation}
\label{claim:two} 
\end{claim}
Observe that, since $\Xb\subset\botud\setminus\Vb$ has full measure 
in $\botud\setminus\Vb$, \eqref{e:claim2} implies~\eqref{e:claim1} 
and concludes the proof of Claim~\ref{claim:one} and hence of Lemma~\ref{l:C1950}.
\paragraph{Proof of Claim~\ref{claim:two}.}
Let $\bpoint\in\Xb\setminus\Eb$ and let $b\in\NN$ such that 
$\afn_{\bpoint}$ is $b\,$-projectively adjacent. 
Then $\bpoint\in\Xb\setminus\Eb(b)$. 
Since $\Xb\setminus\Eb(b)$ is perfect, there is a subset of 
$\Xb\setminus\Eb(b)$, which is at most countable, such that, if 
$\bpoint\in\Xb\setminus\Eb(b)$ does not belong to this subset, 
then there exists a sequence of points of 
$\Xb\setminus\Eb(b)$ that converges to $\bpoint$ from the right 
and one that converges to $\bpoint$ from the left. 

Let $\epsilon>0$. Since the sequence 
${\{g_{(j,b)}\}}_{j\in\NN}$ converges to zero uniformly on $\Xb\setminus\Eb(b)$, 
there exists $n_0\in\NN$ such that if $j\geq n_0$ then 
\begin{equation}
\sup\{(\PoissonI{\cfunction{\Vb}})(\dpoint):\dpoint\in\Gamma_{b}(\bpointx)\cap\ball{\bpointx}{1/j}\}
=
g_{(j,b)}(\bpointx)<\epsilon \quad \forall\bpointx\in\Xb\setminus\Eb(b)
\label{e:sup} 
\end{equation}
Let $r<\frac{2^{-10}}{n_0}$. 
We may assume, without loss of generality, that 
there exists $\theta>0$ such that 
\begin{equation}
\racal{\bpoint}{\theta}\subset\Gamma_{\angolo}^{\ast}(\ballrel{\bpoint}{r})
\label{e:shadow2} 
\end{equation}
and that $\theta<{2^{-10}}{r}$.
Let $\theta_0>0$ with $\theta_0<\theta$ and 
$\bpoint{}e^{i\theta_0}\in\Xb\setminus\Eb(b)$. 
Then~\eqref{e:shadow2} implies that 
$$
\bpoint{}e^{i\theta_0}\in\Gamma_b^*(\dpoint) 
\quad 
\text{ for some }
\dpoint\in\ballrel{\bpoint}{r}
$$
hence $\dpoint\in\Gamma_b(\bpoint{}e^{i\theta_0})\cap\afn_{\bpoint}\cap\ball{\bpoint}{r}$. 
Observe that 
$$
|\dpoint-\bpoint{}e^{i\theta_0}|\leq
|\dpoint-\bpoint{}|+|\bpoint-\bpoint{}e^{i\theta_0}|<r+\theta_0
<\frac{2^{-10}}{n_0}+\frac{2^{-20}}{n_0}
<\frac{1}{n_0}
$$
Since $\dpoint\in\ball{\bpoint{}e^{i\theta_0}}{\frac{1}{n_0}}
\cap{}\Gamma_b(\bpoint{}e^{i\theta_0})$ and 
$\bpoint{}e^{i\theta_0}\in\Xb\setminus\Eb(b)$ 
then~\eqref{e:sup} (with $\bpointx=\bpoint{}e^{i\theta_0}$)
implies 
that  
$$
(\PoissonI{\cfunction{\Vb}})(\dpoint)\leq{}g_{(n_0,b)}(u)<\epsilon
$$
where $\dpoint\in\afn_{\bpoint}\cap\ball{\bpoint}{r}$. 
Thus, we have proved that $\forall\epsilon>0$ 
$\exists\, n_0\in\NN$ such that $\forall \, r\in(0,\frac{2^{-10}}{n_0})$ 
there exists $\dpoint\in\afn_{\bpoint}\cap\ball{\bpoint}{r}$
with 
$(\PoissonI{1_{\Vb}})(\dpoint)<\epsilon$, and this implies that 
\begin{equation}
\liminf_{\stackrel{\dpoint\to\bpoint}{\dpoint\in\afn_{\bpoint}}}\, 
(\PoissonI{\cfunction{\Vb}})(\dpoint) =0,\,
\text{i.e., }
\bpoint\in\Qb. 
\label{e:almost}
\end{equation}
Hence we have proved that $\bpoint\in\Qb$ 
for a.e. $\bpoint\in\Xb\setminus\Eb$.\end{proof}

\paragraph{Proof of Theorem~\ref{thm:maintheorem}}

Observe that, for each $j\in\NN$, $\PoissonI{\cfunction{\Vb_j}} \in \hinftyudone$. 
Lemma~\ref{l:bl1} and 
Lemma~\ref{l:C1950}
%The basic lemmata in Section~\ref{s:basiclemmata} 
and the Carleson tent estimate 
imply that 
$$
\text{for a.e. $\bpoint\in\Cb\setminus\Vb_j$},
\,\,
\osc{\PoissonI{\cfunction{\Vb_j}}}{\bpoint} \geq c_0. 
$$
where $c_0$ is the Carleson tent constant. 
Let $s \eqdef{}1+\frac{1+c_0}{c_0}$. Following 
\cite{Zygmund1949}, 
$$
\bfunction=\sum_{j\geq1}s^{-j}\cfunction{\Vb_j}
$$
It follows that 
$$
\PoissonI{\bfunction}=\sum_{j\geq1}s^{-j}\PoissonI{\cfunction{\Vb_j}}
$$
Define $\Wb\eqdef\cap_{j\geq1}\Vb_j$ and observe that 
$|\Wb|=0$.
\begin{claim}
For a.e. $\bpoint\in\Cb\setminus\Wb$, $\osc{\PoissonI{\bfunction}}{\bpoint}>0$.
\label{claim:1} 
\end{claim}

Since the set $\Cb\setminus\Wb$ has full measure in $\botud$, Claim~\ref{claim:1} completes the proof of Theorem~\ref{thm:maintheorem}. 

\paragraph{Proof of Claim~\ref{claim:1}.}

Let $j$ be the smallest integer $n$ such that $\bpoint \notin \Vb_n$. 
Then $\bpoint$ belongs to the open set 
\begin{equation}
\label{1eqqqnn3}
\bigcap^{j-1}_{k=1} \Vb_k.
\end{equation}
For $k=1,2,\ldots, j-1$, the function $\cfunction{\Vb_k}$ is equal to 1 on the set \eqref{1eqqqnn3}; since this set is open, it follows that 
$\osc{\PoissonI{\cfunction{\Vb_k}}}{\bpoint}=0$
for each $k = 1, 2, \ldots, l-1$.  On the other hand, 
$\osc{s^{-j}\PoissonI{\cfunction{\Vb_j}}
}{\bpoint}\geq s^{-j}c_0$
and
$$
\osc{\sum_{k\geq{}j+1}s^{-k}\PoissonI{(\cfunction{\Vb_k})}}{\bpoint}
\leq
\sum_{k\geq{j+1}}s^{-k}
\leq
s^{-j}\frac{1}{s-1}
$$
It follows that 
$$
\osc{\PoissonI{\bfunction}}{\bpoint}\geq
s^{-j}c_0-s^{-j}\frac{1}{s-1}>0
$$
Since the set 
$(\Cb \setminus \Wb) \setminus \Nb$ 
has measure equal to $2\pi$, the proof of Claim~\ref{claim:1} is completed.

\paragraph{Acknowledgement.}
This research has been supported by the Sida-funded UR-Sweden Program for Research, Higher Learning
and Institution Advancement, sub-program Strengthening Research Capacity in Mathematics, Statistics
and Their Applications.
Partial support from Fondi di Ricerca di Ateneo 
of the Università “G. d'Annunzio” in Chieti and Pescara, 
and from 
Indam-Gnampa, in Italy,  is gratefully acknowledged.

\end{document}